\documentclass[10pt]{report}
\usepackage{amssymb}
\usepackage{latexsym}
\usepackage[dvips]{graphicx}
\usepackage{natbib}
\usepackage{amsmath}
\usepackage{float}
\usepackage{labelfig}
\usepackage{epsfig,psfig}

\def\lemma{\textbf{Lemma}}
\def\theorem{\textbf{Theorem}}

\def\pf{\textbf{Proof: }}

\def\conjecture{\textbf{Conjecture}}
\def\diamond{\diamondsuit}

\begin{document}

\begin{center}
\textbf{Totally geodesic subgraphs of the pants complex.}

Javier Aramayona, Hugo Parlier, Kenneth J. Shackleton

[First draft: August 2006]\\~\\

\end{center}

\setlength{\parindent}{0em}

A\begin{small}BSTRACT\end{small}: Our main theorem asserts that every Farey graph embedded in the $1$-skeleton of the pants complex of any finite type surface is totally geodesic.\\

K\begin{small}EYWORDS\end{small}: Pants complex; Weil-Petersson metric; Farey graph\\

2000 MSC: 57M50 (primary); 05C12 (secondary)\\\\


\setlength{\parindent}{0em}

\textbf{$\S1$. Introduction.\\}

Let $\Sigma$ be a compact, connected and orientable surface, possibly with non-empty boundary, of genus $g(\Sigma)$ and $|\partial \Sigma|$ boundary components, and refer to as the \textit{mapping class group} the group of all self-homeomorphisms of $\Sigma$ up to homotopy.

\setlength{\parindent}{2em}

After Hatcher-Thurston [HT], to the surface $\Sigma$ one may associate a simplicial graph $\mathcal{P}(\Sigma)$, the \textit{pants graph}, whose vertices are all the pants decompositions of $\Sigma$ and any two vertices are connected by an edge if and only if they differ by an elementary move; see $\S2.2$ for an expanded definition. This graph is connected, and one may define a path-metric $d$ on $\mathcal{P}(\Sigma)$ by first assigning length $1$ to each edge and then regarding the result as a length space.

The pants graph, with its own geometry, is a fundamental object to study. Brock [B] revealed deep connections with hyperbolic $3$-manifolds and proved the pants graph is the correct combinatorial model for the Weil-Petersson metric on Teichm\"uller space; the two are quasi-isometric. The isometry group of $(\mathcal{P}, d)$ is also correct in so far as the study of surface groups is concerned, for Margalit [Mar] proved it is almost always isomorphic to the mapping class group of $\Sigma$. In addition, Masur-Schleimer [MasS] proved the pants graph to be one-ended for closed surfaces of genus at least $3$. With only a few exceptions, the pants graph is not hyperbolic in the sense of Gromov [BF].

Our main result concerns the geometry of the pants graph.\\

\newtheorem{1}{\theorem}

\begin{1}
Let $\Sigma$ be a compact, connected and orientable surface, and denote by $\mathcal{F}$ a Farey graph. Let $\phi : \mathcal{F} \longrightarrow \mathcal{P}(\Sigma)$ be a simplicial embedding. Then, $\phi(\mathcal{F})$ is totally geodesic in $\mathcal{P}(\Sigma)$.\\
\label{main}
\end{1}

The completion of the Weil-Petersson metric can be characterised by attaching so-called strata [Mas]. These are totally geodesic subspaces of the completion, by a result of Wolpert [W], and correspond to lower dimensional Teichm\"uller spaces, each with their own Weil-Petersson metric. Combining this with Theorem 1.1 of Brock [B], one finds the Farey subgraphs of the pants graph are uniformly quasi-convex.  Even so, Theorem \ref{main} is not implied by any known coarse geometric result. Moreover, Theorem \ref{main} establishes a complete analogy between the geometry of the Farey subgraphs in a pants graph and the geometry of the corresponding strata lying in the completed Weil-Petersson space.

In order to prove Theorem \ref{main}, we shall need to project paths in the pants graph to paths in the given Farey graph of no greater length. All the notation of Theorem \ref{project} will be explained in $\S2$, but for now we point out the finite set of curves $\pi_{Q}(\nu)$ is the subsurface projection after Masur-Minsky [MasMin] of a pants decomposition $\nu$ to the Farey graph determined by the codimension $1$ multicurve $Q$. The intrinsic metric on this Farey graph, assigning length $1$ to each edge, is denoted by $d_{Q}$.\\

\newtheorem{2}[1]{\theorem}

\begin{2}
Let $\Sigma$ be a compact, connected and orientable surface and denote by $Q$ a codimension $1$ multicurve on $\Sigma$. Let $(\nu_{0}, \ldots, \nu_{n})$ be a path in the pants graph $\mathcal{P}(\Sigma)$. For each index $i \leq n - 1$ and for each $\delta_{i} \in \pi_{Q}(\nu_{i})$, there exists an integer $j \in \{1, 2\}$ and a curve $\delta_{i+j} \in \pi_{Q}(\nu_{i+j})$ such that $d_{Q}(\delta_{i}, \delta_{i+j}) \leq j$.\\
\label{project}
\end{2}

To the authors' knowledge, it has yet to be decided whether there exists a distance non-increasing projection from the whole pants graph to any one of its Farey subgraphs. In the absence of an affirmative result, Theorem \ref{project} may well hold independent interest.

The plan of this paper is as follows. In $\S2$ we recall all the terminology we need, much of which is already standard. In $\S3$ we give an elementary proof to Theorem \ref{project}. Indeed, if $Q$ borders a $1$-holed torus on $\Sigma$, it transpires that one may always take $j = 1$. In $\S4$ we apply Theorem \ref{project} to give an elementary proof to Theorem \ref{main}.

Let us close the introduction by stating the following guiding conjecture.\\

\newtheorem{generalcase}[1]{\conjecture}

\begin{generalcase}

Let $\Sigma_{1}$ and $\Sigma_{2}$ be a pair of compact and orientable surfaces. Let $\phi : \mathcal{P}(\Sigma_{1}) \longrightarrow \mathcal{P}(\Sigma_{2})$ be a simplicial embedding. Then, $\phi(\mathcal{P}(\Sigma_{1}))$ is totally geodesic in $\mathcal{P}(\Sigma_{2})$.\\
\label{conjecture}
\end{generalcase}

\setlength{\parindent}{0em}

A\begin{small}CKNOWLEDGEMENTS:\end{small} The authors wish to thank Brian Bowditch and Koji Fujiwara for their interest and enthusiasm, and to thank Saul Schleimer and Scott Wolpert for their suggestions on the first draft. The authors also wish to thank Centre Interfacultaire Bernoulli EPFL, where the first steps towards proving Theorem \ref{main} were taken, for its hospitality and financial support. The third author was partially supported by a short-term Japan Society for the Promotion of Science postdoctoral fellowship, number P05043, and gratefully acknowledges the JSPS for its financial support. The third author also gratefully acknowledges Institute des Hautes \'Etudes Scientifiques for its hospitality and financial support.

\newpage

\textbf{$\S2$. Background and definitions.}\\

We supply all the background and terminology needed both to understand the statements of our main results, and to make sense of their proofs. Throughout, we regard a loop on $\Sigma$ as the homeomorphic image of a standard circle.\\

\setlength{\parindent}{0em}

\textbf{$\S2.1.$ Curves and multicurves.} A loop on $\Sigma$ is said to be \textit{trivial} if it bounds a disc and \textit{peripheral} if it bounds an annulus whose other boundary component belongs to $\partial \Sigma$. For a non-trivial and non-peripheral loop $c$, we shall denote by $[c]$ its free homotopy class. A \textit{curve} is by definition the free homotopy class of a non-trivial and non-peripheral loop. Given any two curves $\alpha$ and $\beta$, their intersection number $\iota(\alpha, \beta)$ is defined equal to $\rm{min}\{|a \cap b| : a \in \alpha, b \in \beta\}$.

\setlength{\parindent}{2em}

We shall say two curves are \textit{disjoint} if they have zero intersection number, and otherwise say they \textit{intersect essentially}. A pair of curves $\{\alpha, \beta\}$ is said to \textit{fill} the surface $\Sigma$ only if $\iota(\delta, \alpha) + \iota(\delta, \beta) > 0$ for every curve $\delta$. In other words, every curve on $\Sigma$ intersects at least one of $\alpha$ and $\beta$ essentially.

A \textit{multicurve} is a collection of distinct and disjoint curves, and the intersection number for a pair of multicurves is to be defined additively. We denote by $\kappa(\Sigma)$ the size of any maximal multicurve on $\Sigma$, equal to $3g(\Sigma) + |\partial \Sigma| - 3$, and refer to this as the \textit{complexity of $\Sigma$}. Note, the only surfaces of complexity $1$ are the $4$-holed sphere and the $1$-holed torus.

Given a set of disjoint loops $L$, such as the boundary of some subsurface of $\Sigma$, we denote by $[L]$ the multicurve maximal among all multicurves whose every curve is represented by some element of $L$. We shall say a multicurve $\omega$ has \textit{codimension $k$}, for some non-negative integer $k$, only if $|\omega| = \kappa(\Sigma) - k$.\\\\

\begin{figure}[h]
\begin{center}
\AffixLabels{\centerline{\epsfig{file = 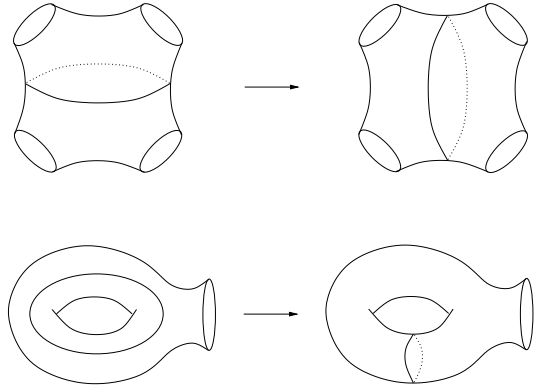, width=7cm, angle= 0}}}
\end{center}
\caption{The two types of elementary move.}
\label{fig:elementarymoves}
\end{figure}

\newpage

\setlength{\parindent}{0em}

\textbf{$\S2.2.$ Pants decompositions.} A \textit{pants decomposition} of a surface is a maximal collection of distinct and disjoint curves, in other words a maximal multicurve. Two pants decompositions $\mu$ and $\nu$ are said to be related by an \textit{elementary move} only if $\mu \cap \nu$ is a codimension $1$ multicurve and the remaining two curves together either fill a $4$-holed sphere and intersect twice or fill a $1$-holed torus and intersect once; consider Figure \ref{fig:elementarymoves} above.\\

\setlength{\parindent}{0em}

\textbf{$\S2.3.$ Arcs.} An \textit{arc} on $\Sigma$ is the homotopy class, relative to $\partial \Sigma$, of an embedded interval ending on $\partial \Sigma$ that does not bound a disc with $\partial \Sigma$. There are broadly two types of arc: those that end on only one component of $\partial \Sigma$, referred to as \textit{waves}, and those that end on two different components of $\partial \Sigma$, referred to as \textit{seams}; see Figure \ref{fig:waveandseam} below.\\\\

\setlength{\parindent}{2em}

\begin{figure}[h]
\begin{center}
\AffixLabels{\centerline{\epsfig{file = 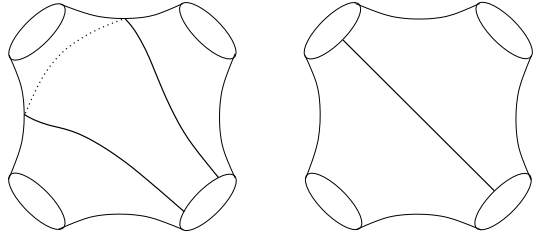, width=7cm, angle= 0}}}
\end{center}
\caption{The two types of arc, respectively a wave and a seam, on a $4$-holed sphere.}
\label{fig:waveandseam}
\end{figure}

Typically, our arcs will live on proper subsurfaces of complexity $1$, noting every arc on a $1$-holed torus is a wave. We may similarly define the intersection number of a pair of arcs, or an arc and a curve, and say two arcs are disjoint or intersect essentially.\\

\setlength{\parindent}{0em}

\textbf{$\S2.4.$ Graphs and paths.} For us, a \textit{path} in a graph shall be a finite sequence of vertices such that any consecutive pair spans an edge; one can recover a topological path by joining up the dots. A \textit{geodesic} is then a distance realising path. Finally, a subgraph $F$ of a metric graph $G$ is said to be \textit{totally geodesic} only if every geodesic in $G$ whose two endpoints belong to $F$ is in fact entirely contained in $F$.\\

\textbf{$\S2.5.$ Farey graphs.} There are numerous ways to build a Farey graph $\mathcal{F}$, any two producing isomorphic graphs. One can start with the rational projective line $\widehat{\Bbb{Q}} := \Bbb{Q} \cup \{\infty\}$, identifying $0$ with $\frac{0}{1}$ and $\infty$ with $\frac{1}{0}$, and take this to be the vertex set of $\mathcal{F}$. Then, two projective rational numbers $\frac{p}{q}, \frac{r}{s} \in \widehat{\Bbb{Q}}$, where $p$ and $q$ are coprime and $r$ and $s$ are coprime, are deemed to span an edge, or $1$-simplex, if and only if $|ps - rq| = 1$. The result is a connected graph in which every edge separates. The graph $\mathcal{F}$ can be naturally represented on a disc; see Figure \ref{fig:farey} below. We shall say a graph \textit{is a Farey graph} if it is isomorphic to $\mathcal{F}$.\\\\

\begin{figure}[h]
\begin{center}
\AffixLabels{\centerline{\epsfig{file = 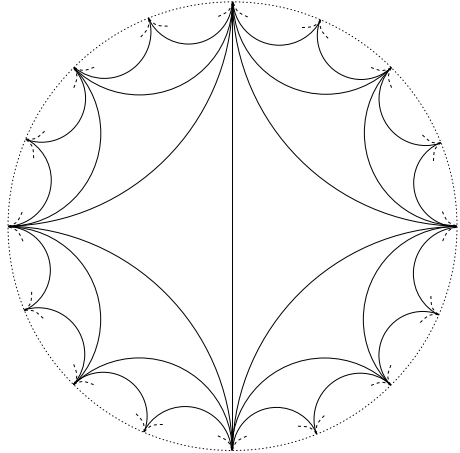, width=6cm, angle= 0}}}
\end{center}
\caption{The Farey graph can be represented on a disc.}
\label{fig:farey}
\end{figure}

\setlength{\parindent}{2em}

It should be noted that both the pants graph of the $4$-holed sphere and the pants graph of the $1$-holed torus are Farey graphs. It follows that any codimension $1$ multicurve $Q$ on $\Sigma$ determines a unique Farey graph $\mathcal{F}_{Q}$ in $\mathcal{P}(\Sigma)$; the converse is Lemma \ref{fareyisco1} from $\S3$. We shall always denote by $d_{Q}$ the intrinsic combinatorial metric on $\mathcal{F}_{Q}$.\\

\setlength{\parindent}{0em}

\textbf{$\S2.6.$ Subsurface projections.} Given a curve $\alpha$ and an incompressible subsurface $Y$ of $\Sigma$, we shall write $\alpha \subset Y$ only if $\alpha$ can be represented by a non-peripheral loop on $Y$. If every loop representing $\alpha$ has non-empty intersection with $Y$ we can say \textit{$\alpha$ and $Y$ intersect}, otherwise we say they are \textit{disjoint}. If every loop representing $\alpha$ intersects $Y$ in at least one interval, we can say \textit{$\alpha$ crosses $Y$}.

\setlength{\parindent}{2em}

For a codimension $1$ multicurve $Q$, let $Y$ denote any complexity $1$ incompressible subsurface of $\Sigma$ such that each curve in $Q$ is disjoint from $Y$. Note then, $Y$ is well defined up to isotopy. Let $\alpha$ be any curve intersecting $Y$, and choose any simple representative $c \in \alpha$ such that $\#(c \cap \partial Y)$ is minimal. We refer to each component of $c \cap Y$ as a \textit{footprint of $c$ on $Y$}, and to the homotopy class of such a footprint as a \textit{footprint of $\alpha$ on $Y$}. Note, footprints of a curve can be arcs or curves.

\begin{figure}[h]
\begin{center}
\AffixLabels{\centerline{\epsfig{file = 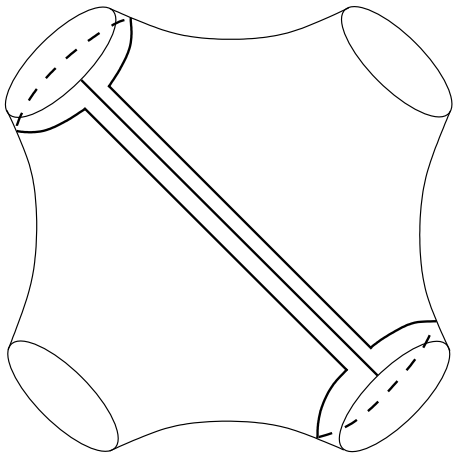, width=6cm, angle= 0}}}
\end{center}
\caption{A seam and a projected curve.}
\label{fig:footprint}
\end{figure}

Given a footprint $b$ for the curve $\alpha$ there only ever exists one curve on $Y$ disjoint from $b$, and such a curve shall be referred to as a \textit{projection of $\alpha$}. Note the set of $\alpha$ projections, each counted once, depends only on $\alpha$ and the original multicurve $Q$, and we denote this set by $\pi_{Q}(\alpha)$. For a second multicurve $\nu$ we may similarly define $\pi_{Q}(\nu)$. The set $\pi_{Q}(\nu)$ is an example of a \textit{subsurface projection}, as defined by Masur-Minsky in $\S1.1$ of [MasMin]. See Figure \ref{fig:footprint} below for an illustration.\\\\

By way of example, we note that $\pi_{Q}(Q) = \emptyset$ and that $\pi_{Q}(\delta) = \{\delta\}$ for any curve $\delta \subset Y$. Finally, if $\delta \subset Y$ is a curve and $\alpha$ is a second curve crossing $Y$ and disjoint from $\delta$, then $\delta \in \pi_{Q}(\alpha)$.\\\\

\setlength{\parindent}{0em}

\textbf{$\S3$. One proof of Theorem \ref{project}.}\\

Let us start with two well known results, the first characterising the Farey subgraphs of any given pants graph and the second relating low intersection numbers to distances for a pair of curves on a $4$-holed sphere.\\

\newtheorem{3}[1]{\lemma}

\begin{3}
Let $\phi : \mathcal{F} \longrightarrow \mathcal{P}(\Sigma)$ be a simplicial embedding. Then, there exists a codimension $1$ multicurve on $\Sigma$ contained in every vertex of $\phi(\mathcal{F})$.\label{fareyisco1}
\end{3}

\pf This follows from the following two remarks. First, note the vertices of any $3$-cycle from $\phi(\mathcal{F})$ always intersect in a common codimension $1$ multicurve. Second, note for any two vertices $\mu$ and $\nu$ of $\phi(\mathcal{F})$ there exists a finite sequence of $3$-cycles $\Delta_{0}, \ldots, \Delta_{n}$ such that $\mu$ is a vertex of $\Delta_{0}$, such that $\nu$ is a vertex of $\Delta_{n}$, and such that $\Delta_{i} \cap \Delta_{i+1}$ is an edge for each index $i$. One can then prove Lemma \ref{fareyisco1} by an induction. $\diamond$\\

\newtheorem{3+}[1]{\lemma}

\begin{3+}
Let $Y$ be a $4$-holed sphere. Any two vertices $\delta_{0}, \delta_{2}$ of $\mathcal{P}(Y)$ of intersection number at most $4$ are at distance $d(\delta_{0}, \delta_{2})$ at most $2$.
\label{lowintlowdist}
\end{3+}

\pf There exists a curve $\delta_{1}$ on $Y$ such that $\iota(\delta_{1}, \delta_{j}) \leq 2$ for both $j \in \{0, 2\}$; if $\iota(\delta_{0}, \delta_{2}) = 4$ then such a curve can be explicitly constructed by performing an elementary surgery on either of $\delta_{0}$ or $\delta_{2}$. It follows that $d(\delta_{0}, \delta_{2}) \leq d(\delta_{0}, \delta_{1}) + d(\delta_{1}, \delta_{2}) \leq 1 + 1 = 2$. $\diamond$\\

\setlength{\parindent}{2em}

The following two lemmata shall be applied in what will become known as Case B for the $4$-holed sphere, Lemma \ref{3curve} playing an especially important r\^ole.\\

\newtheorem{4}[1]{\lemma}
\begin{4}
Let $P$ be a pants decomposition of $\Sigma$, and let $Y$ be a connected complexity $1$ incompressible subsurface of $\Sigma$. If $P$ does not contain $[\partial Y]$, then $P$ contains at least two curves intersecting $Y$.\label{2curve}
\end{4}

\setlength{\parindent}{0em}

\pf We shall denote by $\kappa_{*}(Y)$ the size of a maximal multicurve on $\Sigma$ whose every curve does not intersect $Y$. Let $\omega \subset P$ be the set of all curves in $P$ that do not intersect $Y$. We have $$|P| = \kappa(\Sigma) = \kappa(Y) + \kappa_{*}(Y) = 1 + \kappa_{*}(Y) \geq 1 + |\omega| + 1 = 2 + |\omega|,$$ and so $|P| \geq 2 + |\omega|$ as required. $\diamond$\\

\newtheorem{5}[1]{\lemma}

\begin{5}
Let $P$ be a pants decomposition of $\Sigma$, and let $Y$ be an incompressible subsurface of $\Sigma$ homeomorphic to a $4$-holed sphere. If there exist two distinct curves in $[\partial Y]$ not contained in $P$, then $P$ contains at least three curves intersecting $Y$.\label{3curve}
\end{5}

\setlength{\parindent}{0em}

\pf Let $\omega \subset P$ be the set of all curves in $P$ that do not intersect $Y$. We have $$|P| = \kappa(\Sigma) = \kappa(Y) + \kappa_{*}(Y) = 1 + \kappa_{*}(Y) \geq 1 + |\omega| + 2 = 3 + |\omega|,$$ and so $|P| \geq 3 + |\omega|$ as required. $\diamond$\\

\setlength{\parindent}{2em}

We now turn to proving Theorem \ref{project}, denoting by $Y$ the complexity $1$ subsurface of $\Sigma$ complementary to $Q$. Note the statement of Theorem \ref{project} holds vacuously if $\kappa(\Sigma) \leq 0$ and trivially if $\kappa(\Sigma) = 1$, since then $\phi$ is an isomorphism. When $\kappa(\Sigma) = 2$, the surface $\Sigma$ is either a $5$-holed sphere or a $2$-holed torus. If the genera $g(Y)$ and $g(\Sigma)$ are equal, then each footprint of $\nu_{i+1}$ on $Y$ is therefore either a curve or a wave. As such, there exists a curve $\delta_{i+1} \in \pi_{Q}(\nu_{i+1})$ such that $\delta_{i}$ and $\delta_{i+1}$ are either equal or intersect minimally. We can then take $j = 1$, noting $d_{Q}(\delta_{i}, \delta_{i+1}) = 1$. The remaining case, $\Sigma$ the $2$-holed torus and $Y$ the $4$-holed sphere, is deferred to Appendix.

For the remainder of this section, it is to be assumed that $\kappa(\Sigma) \geq 3$. Let $\delta_{i} \in \pi_{Q}(\nu_{i})$. In constructing a curve $\delta_{i+1}$ or $\delta_{i+2}$, as per the statement of Theorem \ref{project}, we note Lemma \ref{fareyisco1} tells us it is enough to consider separately the case $Y$ is a $4$-holed sphere and the case $Y$ is a $1$-holed torus.\\


\setlength{\parindent}{0em}

\begin{center}
\begin{small}
$Y$ IS A $4$-HOLED SPHERE\\
\end{small}
\end{center}

The case $Y$ is a $4$-holed sphere separates into two main cases, according as $\delta_{i}$ belongs to $\nu_{i}$ or does not belong to $\nu_{i}$.\\

\textbf{Case A:} $\delta_{i} \in \nu_{i}$.\\

\textbf{I.} $\delta_{i} \in \nu_{i+1}$. Take $j = 1$ and $\delta_{i+1} = \delta_{i}$.\\

\textbf{II.} $\delta_{i} \notin \nu_{i+1}$. We may still take $j = 1$ and choose any $\delta_{i+1} \in \pi_{Q}(\nu_{i+1})$, for $\delta_{i}$ is a curve and, as such, is disjoint from $[\partial Y]$. Now $\iota(\delta_{i}, \delta_{i+1}) \leq 2$ and so $d_{Q}(\delta_{i}, \delta_{i+1}) \leq 1$.\\

\textbf{Case B:} $\delta_{i} \notin \nu_{i}$.\\

There exists a $Y$-footprint $a_{i}$ of $\nu_{i}$ such that $\iota(a_{i}, \delta_{i}) = 0$, that is $\delta_{i}$ is uniquely determined by $a_{i}$. We denote by $\alpha_{i}$ any curve from $\nu_{i}$ having $a_{i}$ as a footprint. Let $a_{i+1}$ be any footprint of $\nu_{i+1}$ on $Y$, and let $\alpha_{i+1}$ be any element of $\nu_{i+1}$ having $a_{i+1}$ as a footprint.\\

\textbf{I.} \textit{$a_{i}$ and $a_{i+1}$ intersect essentially.} Since $a_{i}$ and $a_{i+1}$ intersect essentially, so must the two curves $\alpha_{i}$ and $\alpha_{i+1}$. Moreover, as $\delta_{i} \notin \nu_{i}$, so $a_{i}$ can only be an arc.

\setlength{\parindent}{2em}

Suppose first that $a_{i+1}$ is a curve. Then, $\alpha_{i+1}$ and $a_{i+1}$ are equal. According to Lemma \ref{2curve} there exists a curve $\alpha_{i}' \in \nu_{i}$ such that $\alpha_{i} \neq \alpha_{i}'$ and such that $\alpha_{i}'$ intersects $Y$. Since $d(\nu_{i}, \nu_{i+1}) = 1$ and since $\iota(\alpha_{i}, \alpha_{i+1}) \neq 0$, so $\alpha_{i}' \in \nu_{i+1}$. The set $\{\alpha_{i}', \alpha_{i+1}\} \cap \nu_{i+2}$ is therefore non-empty. Let $\gamma \in \{\alpha_{i}', \alpha_{i+1}\} \cap \nu_{i+2}$ and take $j = 2$. There exists $\delta_{i+2} \in \pi_{Q}(\gamma)$ such that $\iota(\delta_{i}, \delta_{i+2}) \leq 4$ and so, according to Lemma \ref{lowintlowdist}, $d_{Q}(\delta_{i}, \delta_{i+2}) \leq 2$.

\setlength{\parindent}{2em}

Henceforth, $a_{i+1}$ shall always be an arc. Appealing to Lemma \ref{2curve}, there exists a $Y$-footprint $a_{i+1}'$ of $\nu_{i+1}$ and a corresponding curve $\alpha_{i+1}' \in \nu_{i+1}$ such that $\alpha_{i+1}' \neq \alpha_{i+1}$. Since $d(\nu_{i}, \nu_{i+1}) = 1$ it follows that $\iota(a_{i}, a_{i+1}') = 0$. Note, if $a_{i+1}'$ is a curve then $\alpha_{i+1}' = a_{i+1}'$ and we may take $j = 1$ and $\delta_{i+1} = \alpha_{i+1}'$.

Henceforth, $a_{i+1}'$ is assumed to be an arc. We observe $a_{i+1}$ and $a_{i+1}'$ are distinct arcs, since $a_{i+1}$ intersects $a_{i}$ essentially whereas $a_{i+1}'$ is disjoint from $a_{i}$. The first case, B.I., will now be completed by considering in turn the two topological possibilities for $a_{i}$.\\

\setlength{\parindent}{0em}

(i) \textit{$a_{i}$ is a wave.} Let $\gamma_{i+1}' \in \pi_{Q}(\nu_{i+1})$ be such that $\iota(\gamma_{i+1}', a_{i+1}') = 0$. Then, $\iota(\delta_{i}, \gamma_{i+1}') \leq 2$ and therefore $d_{Q}(\delta_{i}, \gamma_{i+1}') \leq 1$. We can thus take $j = 1$ and $\delta_{i+1} = \gamma_{i+1}'$.\\

(ii) \textit{$a_{i}$ is a seam.} Let us suppose the two components of $\partial Y$ on which $a_{i}$ ends are not homotopic on $\Sigma$; the remaining case seems to require individual consideration, and so we prefer to postpone this to Appendix.

\setlength{\parindent}{2em}

Suppose $\{a_{i}, a_{i+1}'\}$ ends on at least three different components of $\partial Y$. Let $\gamma_{i+1}' \in \pi_{Q}(\nu_{i+1})$ be such that $\iota(\gamma_{i}', a_{i+1}') = 0$. Then, $\iota(\delta_{i}, \gamma_{i+1}') \leq 2$ and therefore $d_{Q}(\delta_{i}, \gamma_{i+1}') \leq 1$. We now take $j = 1$ and $\delta_{i+1} = \gamma_{i+1}'$.

Suppose instead $\{a_{i}, a_{i+1}'\}$ now ends on at most two, and therefore exactly two, different components of $\partial Y$. Since $a_{i+1}'$ is a seam and since $\alpha_{i+1} \in \nu_{i+1}$, so $\nu_{i+1}$ fails to contain at least two curves from $[\partial Y]$. Appealing to Lemma \ref{3curve}, there exists a curve $\alpha_{i+1}'' \in \nu_{i+1}$ such that $\alpha_{i+1}'' \notin \{\alpha_{i+1}, \alpha_{i+1}'\}$ and such that $\alpha_{i+1}''$ intersects $Y$. Since $d(\nu_{i}, \nu_{i+1}) = 1$ and since $\iota(\alpha_{i}, \alpha_{i+1}) \neq 0$, so $\iota(\alpha_{i}, \alpha_{i+1}'') = 0$. Moreover, since $\alpha_{i} \notin \nu_{i+1}$, so $\alpha_{i} \neq \alpha_{i+1}''$. As $d(\nu_{i+1}, \nu_{i+2}) = 1$, so $\{\alpha_{i+1}', \alpha_{i+1}''\} \cap \nu_{i+2} \neq \emptyset$. Let $\gamma \in \{\alpha_{i+1}', \alpha_{i+1}''\} \cap \nu_{i+2}$. We now take $j = 2$ and $\delta_{i+2} \in \pi_{Q}(\gamma)$, noting that $\iota(\delta_{i}, \delta_{i+2}) \leq 4$ and, as such, $d_{Q}(\delta_{i}, \delta_{i+2}) \leq 2$.\\

\setlength{\parindent}{0em}

\textbf{II.} \textit{$a_{i}$ and $a_{i+1}$ are disjoint.} First note that, if either of $a_{i}$ and $a_{i+1}$ is a wave, we may take $j = 1$ and $\delta_{i+1} \in \pi_{Q}(\alpha_{i+1})$ such that $\iota(\delta_{i+1}, a_{i+1}) = 0$. Then, $\iota(\delta_{i}, \delta_{i+1}) \leq 2$ and, as such, $d_{Q}(\delta_{i}, \delta_{i+1}) \leq 1$. Henceforth, we assume that $a_{i}$ and $a_{i+1}$ are both seams.

\setlength{\parindent}{2em}

If $\{a_{i}, a_{i+1}\}$ ends on at least three components of $\partial Y$ we may take $j = 1$ and $\delta_{i+1} \in \pi_{Q}(\alpha_{i+1})$ such that $\iota(\delta_{i+1}, a_{i+1}) = 0$. Then, $\iota(\delta_{i}, \delta_{i+1}) \leq 2$ and, as such, $d_{Q}(\delta_{i}, \delta_{i+1}) \leq 1$.

Thus, we may assume that $\{a_{i}, a_{i+1}\}$ ends on at most two, and therefore exactly two, components of $\partial Y$. By assumption, $\nu_{i+1}$ does not contain $[\partial Y]$. According to Lemma \ref{2curve}, there exists a second $Y$-footprint $a_{i+1}'$ for some curve $\alpha_{i+1}' \in \nu_{i+1}$ such that $\alpha_{i+1}$ and $\alpha_{i+1}'$ are distinct. If $a_{i}$ and $a_{i+1}$ are equal then $\delta_{i} \in \pi_{Q}(\nu_{i+1})$, and we may take $j = 1$ and $\delta_{i+1} = \delta_{i}$. We may therefore assume $a_{i}$ and $a_{i+1}'$ are not equal.

\setlength{\parindent}{2em}

If $a_{i}$ and $a_{i+1}'$ intersect essentially, then we may appeal to Case B.I. with $a_{i+1}'$ substituted for $a_{i+1}$. We may thus assume that $a_{i}$ and $a_{i+1}'$ are disjoint. Since three disjoint arcs on $Y$ cannot end on at most two components of $\partial Y$, it follows $\{a_{i}, a_{i+1}'\}$ ends on at least three different components of $\partial Y$. We can now take $j = 1$ and $\delta_{i+1} \in \pi_{Q}(\alpha_{i+1}')$ such that $\iota(\delta_{i+1}, a_{i+1}') = 0$.\\

This concludes the case $Y$ is a $4$-holed sphere.\\

\setlength{\parindent}{0em}

\begin{center}
\begin{small}
$Y$ IS A $1$-HOLED TORUS\\
\end{small}
\end{center}

The case of the $1$-holed torus is more straightforward, for here each arc is a wave, and can be treated by considering separately four mutually exclusive cases.\\

\textbf{I.} \textit{$\nu_{i}, \nu_{i+1}$ contain $[\partial Y]$.} Let $\delta_{i+1}$ denote the only curve contained in $\pi_{Q}(\nu_{i+1})$. We may then take $j = 1$ and note $d_{Q}(\delta_{i}, \delta_{i+1})) \leq 1$.\\

\textbf{II.} \textit{$\nu_{i}$ contains $[\partial Y]$, whereas $\nu_{i+1}$ does not.} Then, $\delta_{i} \in \nu_{i+1}$. We may take $j = 1$ and $\delta_{i+1} = \delta_{i}$.\\

\textbf{III.} \textit{$\nu_{i+1}$ contains $[\partial Y]$, whereas $\nu_{i}$ does not.} Then, $\nu_{i+1}$ contains a single curve $\gamma_{i+1}$ such that $\gamma_{i+1} \subset Y$. Since $d(\nu_{i}, \nu_{i+1}) = 1$, so $\gamma_{i+1} \in \nu_{i}$ and hence $\gamma_{i+1} \in \pi_{Q}(\nu_{i})$. As $\pi_{Q}(\nu_{i})$ contains only one element, so $\gamma_{i+1} = \delta_{i}$. We may now take $j = 1$ and $\delta_{i+1} = \delta_{i}$.\\

\textbf{IV.} \textit{Neither $\nu_{i}$ nor $\nu_{i+1}$ contains $[\partial Y]$.} There exists a $Y$-footprint $a_{i}$ of $\nu_{i}$ such that $\iota(\delta_{i}, a_{i}) = 0$. According to Lemma \ref{2curve}, there exist two footprints $a_{i+1}$ and $a_{i+1}'$ of $\nu_{i+1}$ corresponding to different elements of $\nu_{i+1}$. Since $d(\nu_{i}, \nu_{i+1}) = 1$, so at least one of these footprints, say $a_{i+1}$, is disjoint from $a_{i}$. We may take $j = 1$ and $\delta_{i+1} \in \pi_{Q}(\nu_{i+1})$ such that $\iota(\delta_{i+1}, a_{i+1}) = 0$. Note, $\iota(\delta_{i}, \delta_{i+1}) \leq 1$ and, as such, $d_{Q}(\delta_{i}, \delta_{i+1}) \leq 1$.\\

\setlength{\parindent}{2em}

This concludes the case $Y$ is a $1$-holed torus, thus concluding a proof of Theorem \ref{project}. $\diamond$\\\\


\setlength{\parindent}{0em}

\textbf{$\S4$. One proof of Theorem \ref{main}.}\\

Let $\mathcal{F}$ be a Farey graph and let $\phi : \mathcal{F} \longrightarrow \mathcal{P}(\Sigma)$ be a simplicial embedding. There exists a unique codimension $1$ multicurve $Q$ on $\Sigma$ such that $Q$ is contained in every vertex of $\phi(\mathcal{F})$; see Lemma \ref{fareyisco1}.

\setlength{\parindent}{2em}

Suppose, for contradiction, that $\phi(\mathcal{F})$ is not totally geodesic. Then, there exist two vertices $\mu$ and $\nu$ of $\phi(\mathcal{F})$ and a geodesic $\mu = \nu_{0}, \nu_{1}, \ldots, \nu_{n} = \nu$ in $\mathcal{P}(\Sigma)$ not entirely contained in $\phi(\mathcal{F})$. Let $i$ be the minimal index such that $\nu_{i} \notin \phi(\mathcal{F})$, noting $1 \leq i \leq n-1$. Let $\delta_{i-1}$ and $\delta_{i}$ be, respectively, the one element of $\pi_{Q}(\nu_{i-1})$ and the one element of $\pi_{Q}(\nu_{i})$, noting that $\delta_{i-1} = \delta_{i}$. According to Theorem \ref{project} there exists a sequence of integers $(n_{j}) \subseteq \{i-1, \ldots, n\}$, containing $i-1$ and $n$, and a corresponding sequence of curves $\delta_{n_{j}} \in \pi_{Q}(\nu_{n_{j}})$ such that $0 < n_{j+1} - n_{j} \leq 2$, for each $j$, and such that $d_{Q}(\delta_{n_{j}}, \delta_{n_{j+1}}) \leq n_{j+1} - n_{j}$, for each $j$. Necessarily, $\phi(\delta_{i-1}) = \nu_{i-1}$ and $\phi(\delta_{n}) = \nu_{n}$. We note that $$d_{Q}(\delta_{i-1}, \delta_{n}) = d_{Q}(\delta_{i}, \delta_{n}) \leq \sum_{j} d_{Q}(\delta_{n_{j}}, \delta_{n_{j+1}}) \leq \sum_{j} n_{j+1} - n_{j} = n - i,$$ and, since paths in $\mathcal{F}$ determine paths in $\mathcal{P}(\Sigma)$ via $\phi$, so it follows that $$d(\nu_{0}, \nu_{n}) = d(\nu_{0}, \nu_{i-1}) + d(\nu_{i-1}, \nu_{n}) \leq i - 1 + d_{Q}(\delta_{i-1}, \delta_{n}) \leq i - 1 + n - i = n-1.$$ To be more succinct, $d(\nu_{0}, \nu_{n}) \leq n - 1$. This is a contradiction, and the statement of Theorem \ref{main} follows. $\diamond$\\\\

\newpage

\setlength{\parindent}{0em}

\textbf{Appendix.\\}

We treat separately one instance of Case B.I(ii) from the proof of Theorem \ref{project}, where the seam $a_{i}$ ends on two distinct but homotopic components of $\partial Y$. This simultaneously treats the case $\Sigma$ is a $2$-holed torus and $Y$ is a $4$-holed sphere. In either instance, we cannot appeal to Lemma \ref{3curve}. We recall $a_{i+1}$ is a footprint of $\alpha_{i+1} \in \nu_{i+1}$ on $Y$ that intersects $a_{i}$ essentially, and that $a_{i+1}'$ is a footprint of $\alpha_{i+1}' \in \nu_{i+1}$ on $Y$ both disjoint from and non-homotopic to $a_{i}$. In addition, we might as well assume $\{a_{i}, a_{i+1}\}$ and $\{a_{i}, a_{i+1}'\}$ both end on precisely two distinct components of $\partial Y$, for we may otherwise take $j = 1$ and readily find $\delta_{i+1} \in \pi_{Q}(\nu_{i+1})$ as claimed.\\

\begin{figure}[h]
\leavevmode \SetLabels
\L(0.49*1.035) $a_{i}$\\
\L(0.35*1.035) $a_{i+1}$\\
\L(0.27*0.81) $a_{i+1}'$\\
\L(0.68*0.92) $\gamma'$\\
\endSetLabels
\begin{center}
\AffixLabels{\centerline{\epsfig{file = 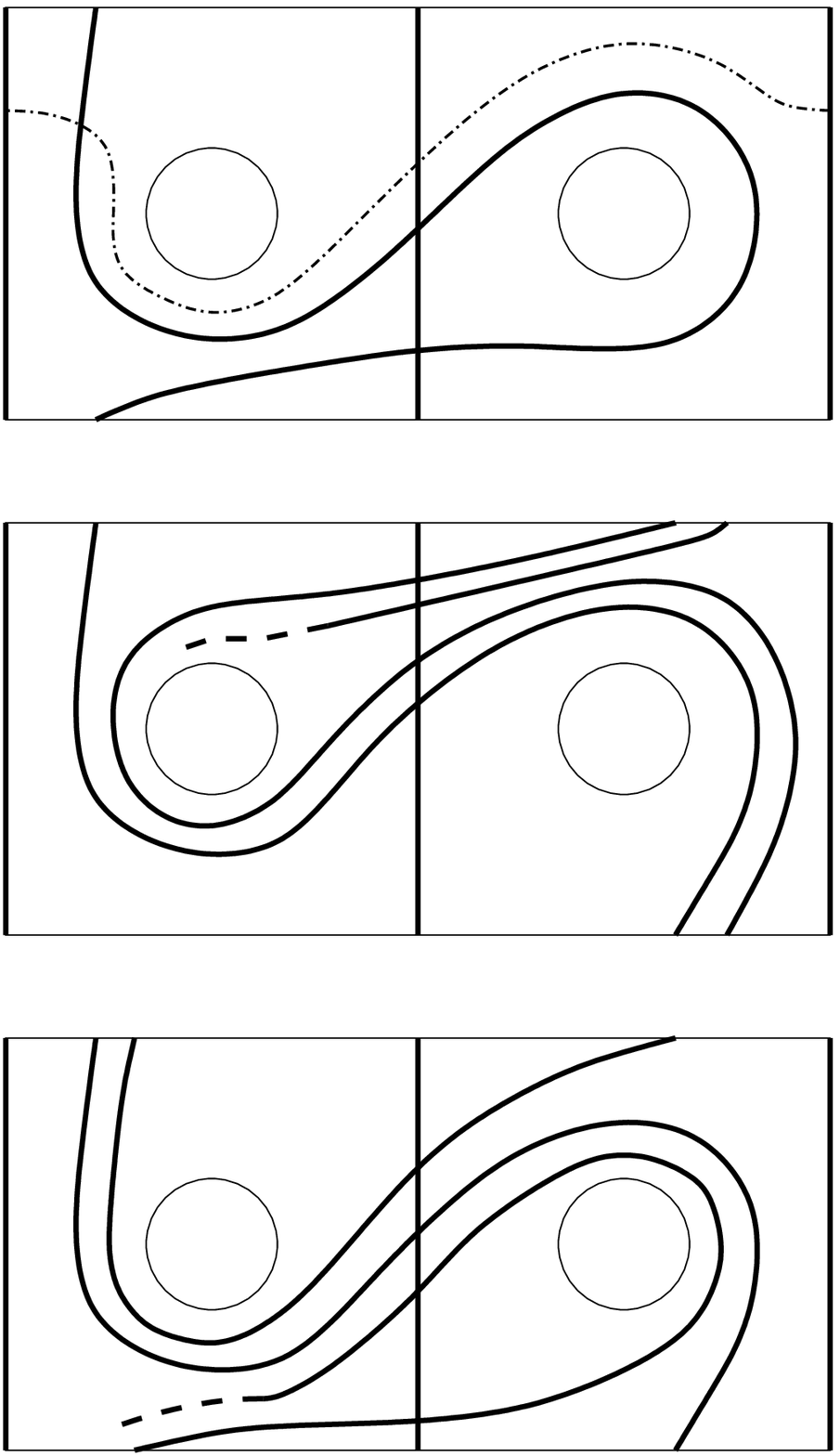, width=4cm, angle= 0}}}
\end{center}
\caption{The case $a_{i+1}$ is a seam, giving three possibilities, up to symmetry, of which only the first can be legal. We represent $Y$ as a disc with four corners and two holes, identifying the left and right vertical edges to give $a_{i+1}'$ and the middle edge with $a_{i}$. The top and bottom edges correspond to distinct components of $\partial Y$, homotopic on $\Sigma$.}
\label{fig:specialcases1}
\end{figure}

\textbf{I.} \textit{$a_{i+1}$ is a seam.} Up to symmetry there are only three possibilities, as per Figure \ref{fig:specialcases1}. In both the second and the third of these, $\alpha_{i}$ and $\alpha_{i+1}$ necessarily have intersection number at least $3$. It follows $\iota(\nu_{i}, \nu_{i+1}) \geq 3$. However, $d(\nu_{i}, \nu_{i+1}) = 1$ and, as such, $\iota(\nu_{i}, \nu_{i+1}) \leq 2$. This is a contradiction.\\\\

\setlength{\parindent}{2em}

The first case is somewhat different to any other treated in this paper. We note $\iota(\alpha_{i}, \alpha_{i+1}) = 2$. Let $\gamma \in \pi_{Q}(\alpha_{i+1})$ be the curve such that $\iota(\gamma, a_{i+1}) = 0$. Then, $\iota(\delta_{i}, \gamma) = 8$. However, there exists a further curve $\gamma' \subset Y$ such that $\iota(\delta_{i}, \gamma') = 2$ and $\iota(\gamma', \gamma) = 2$. It follows $d_{Q}(\delta_{i}, \gamma) \leq 2$, in fact precisely $2$; see Figure \ref{fig:specialcases1}. We may therefore take $j = 2$ and find $\delta_{i+2} \in \pi_{Q}(\{\alpha_{i+1}, \alpha_{i+1}'\} \cap \nu_{i+2})$ such that $d_{Q}(\delta_{i}, \delta_{i+2}) \leq 2$.\\

\setlength{\parindent}{0em}

\begin{figure}[h]
\leavevmode \SetLabels
\L(0.49*1.035) $a_{i}$\\
\L(0.35*1.035) $a_{i+1}$\\
\L(0.27*0.84) $a_{i+1}'$\\
\endSetLabels
\begin{center}
\AffixLabels{\centerline{\epsfig{file = 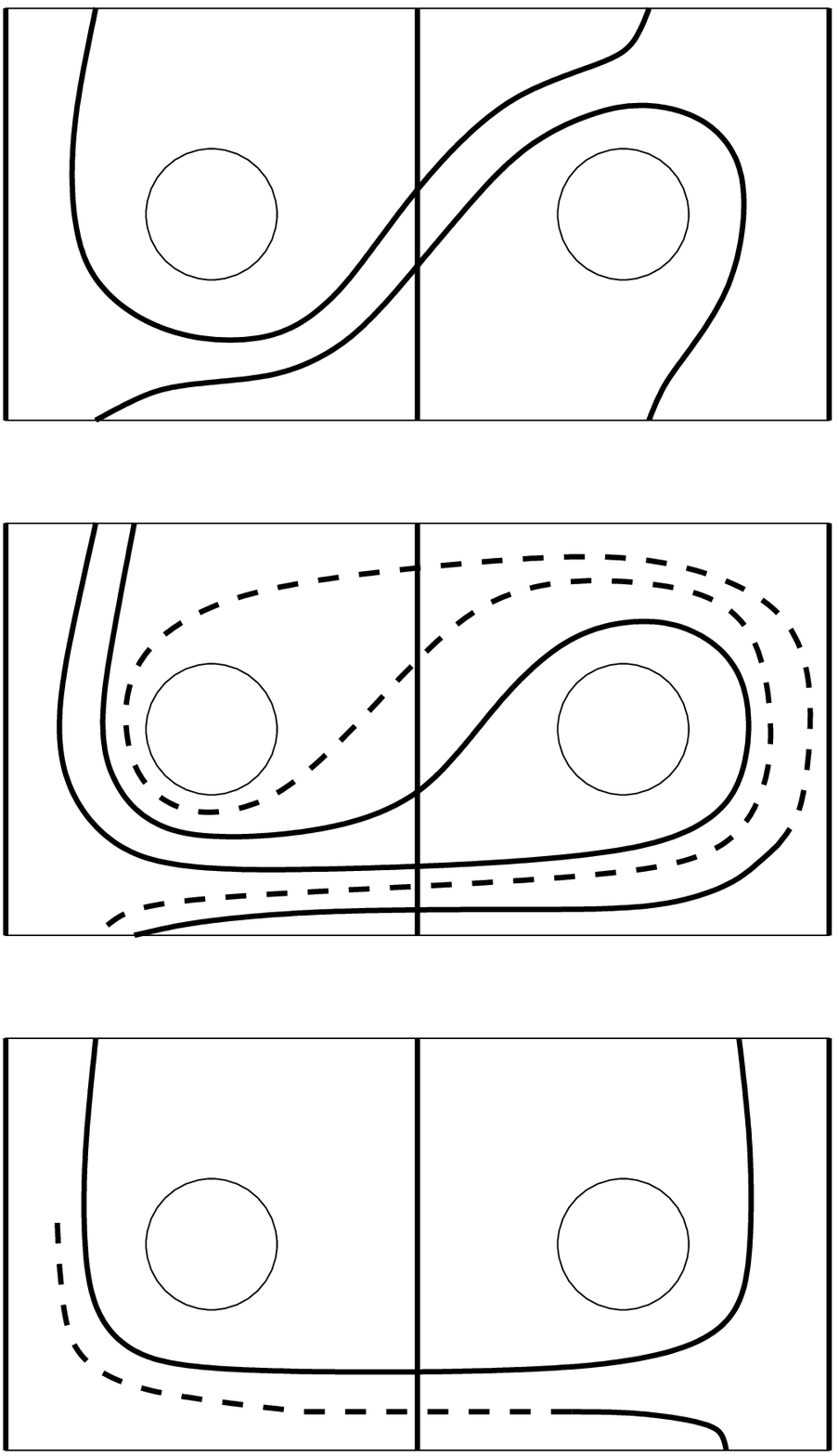, width=4cm, angle= 0}}}
\end{center}
\caption{The case $a_{i+1}$ is a wave, giving three possibilities up to symmetry. In fact, the second and third possibilities are both illegal. For example, in the second diagram we have $\iota(\alpha_{i}, \alpha_{i+1}) \geq 3$.}
\label{fig:specialcases2}
\end{figure}

\textbf{II.} \textit{$a_{i+1}$ is a wave.} We keep Figure \ref{fig:specialcases2} in mind. Let $\gamma \in \pi_{Q}(\nu_{i+1})$ be such that $\iota(\gamma, a_{i+1}) = 0$. Since $d(\nu_{i+1}, \nu_{i+2}) = 1$, so the set $\{\alpha_{i+1}, \alpha_{i+1}'\} \cap \nu_{i+2}$ is non-empty. We take $j = 2$ and let $\delta_{i+2} \in \pi_{Q}(\{\alpha_{i+1}, \alpha_{i+1}'\} \cap \nu_{i+2})$, noting $\iota(\delta_{i}, \delta_{i+2}) \leq 4$ and, as such, $d_{Q}(\delta_{i}, \delta_{i+2}) \leq 2$.\\\\

\newpage

\textbf{References.}\\

\setlength{\parindent}{0em}
\setlength{\parskip}{0.5em plus 0.5em minus 0.5em}

[B] J. F. Brock, \textit{The Weil-Petersson metric and volumes of $3$-dimensional hyperbolic convex cores} : Journal of the American Mathematical Society \textbf{16} No. 3 (2003) 495--535.

[BF] J. F. Brock, B. Farb, \textit{Curvature and rank of Teichm\"uller space} : American Journal of Mathematics \textbf{128} (2006) 1--22.

[HT] A. E. Hatcher, W. P. Thurston, \textit{A presentation for the mapping class group of a closed orientable surface} : Topology \textbf{19} (1980) 221--237.

[Mar] D. Margalit, \textit{Automorphisms of the pants complex} : Duke Mathematical Journal \textbf{121} No. 3 (2004).

[Mas] H. A. Masur, \textit{Extension of the Weil-Petersson metric to the boundary of Teichmuller space} : Duke Mathematical Journal \textbf{43} no. 3 (1976) 623--635.

[MasMin] H. A. Masur, Y. N. Minsky, \textit{Geometry of the complex of curves II: Hierarchical structure} : Geometry \& Functional Analysis \textbf{10} (2000) 902--974.

[MasS] H. A. Masur, S. Schleimer, \textit{The pants complex has only one end} : to appear in ``Spaces of Kleinian groups'' (eds. Y. N. Minsky, M. Sakuma, C. M. Series), Proceedings of Newton Institute conference, London Mathematical Society Lecture Notes Series.

[W] S. A. Wolpert, \textit{Geometry of the Weil-Petersson completion of Teichm\"uller space} : Surveys in Differential Geometry VIII: Papers in honor of Calabi, Lawson, Siu and Uhlenbeck, editor S. T. Yau, International Press (2003).\\

\setlength{\parskip}{0em plus 0.5em minus 0.5em}

\begin{small}
Javier Aramayona\\
Mathematics Institute\\
University of Warwick\\
Coventry CV4 7AL\\
England\\
homepage: http://www.maths.warwick.ac.uk/$\sim$jaram\\
e-mail: jaram@maths.warwick.ac.uk\\

Hugo Parlier\\
Section de Math\'ematiques\\
Universit\'e de Gen\`eve\\
1211 Gen\`eve 4\\
Suisse\\
homepage: http://www.unige.ch/math/folks/parlier/\\
e-mail: hugo.parlier@math.unige.ch\\

Kenneth J. Shackleton (corresponding author)\\
Department of Mathematical and Computing Sciences\\
Tokyo Institute of Technology\\
2-12-1 O-okayama\\
Meguro-ku\\
Tokyo 152-8552\\
Japan\\
homepage: http://www.maths.soton.ac.uk/$\sim$kjs\\
e-mail: kjs2006@alumni.soton.ac.uk
\end{small}

\end{document}